\documentclass[%
 reprint,
 amsmath,amssymb,
 aps,
]{revtex4-1}

\usepackage{graphicx}
\usepackage{dcolumn}
\usepackage{bm}

\usepackage[T1]{fontenc}
\usepackage[utf8]{inputenc}

\usepackage{amssymb,amsfonts,amsmath,mathtools,setspace,savesym,natbib}

\savesymbol{iint}
\savesymbol{iiint}

\newcommand\ket[1]{\left|#1\right>}
\newcommand\vev[1]{\left<#1\right>}

\newcommand{\IR}{\mathbb{R}}

\newcommand*\xbar[1]{
\hbox{
\vbox{
\hrule height 0.66pt
\kern0.33ex
\hbox{
\kern-0.2em
\ensuremath{#1}
\kern-0.0em
}
}
}
}

\newcommand*\blankbar[1]{
\hbox{
\vbox{
\hrule height 0.0pt
\kern0.66ex
\hbox{
\kern-0.2em
\ensuremath{#1}
\kern-0.0em
}
}
}
}

\begin{document}

\preprint{APS/123-QED}

\title{A quantum framework for likelihood ratios}

\author{Rachael L. Bond}
\email{rachael.be@rachael.be}
\affiliation{The School of Psychology, The University of Sussex, Falmer, East Sussex, BN1 9QH, United Kingdom.}

\author{Yang-Hui He}
\email{hey@maths.ox.ac.uk}
\affiliation{Department of Mathematics, City University, London, EC1V 0HB, United Kingdom, and}
 \affiliation{School of Physics, NanKai University, Tianjin, 300071, P.R. China, and}
 \affiliation{Merton College, Oxford University, Oxford, OX1 4JD, United Kingdom.}

\author{Thomas C. Ormerod}
\email{t.ormerod@sussex.ac.uk}
\affiliation{The School of Psychology, The University of Sussex, Falmer, East Sussex, BN1 9QH, United Kingdom.}

\date{\today}

\begin{abstract}
The ability to calculate precise likelihood ratios is fundamental to many STEM areas, such as decision-making theory, biomedical science, and engineering. However, there is no assumption-free statistical methodology to achieve this. For instance, in the absence of data relating to covariate overlap, the widely used Bayes' theorem either defaults to the marginal probability driven ``naive Bayes' classifier'', or requires the use of compensatory expectation-maximization techniques. Equally, the use of alternative statistical approaches, such as multivariate logistic regression, may be confounded by other axiomatic conditions, e.g., low levels of co-linearity. This article takes an information-theoretic approach in developing a new statistical formula for the calculation of likelihood ratios based on the principles of quantum entanglement. In doing so, it is argued that this quantum approach demonstrates: that the likelihood ratio is a real quality of statistical systems; that the naive Bayes' classifier is a special case of a more general quantum mechanical expression; and that only a quantum mechanical approach can overcome the axiomatic limitations of classical statistics. 
\end{abstract}

\pacs{89}
\maketitle

In recent years, Bayesian statistical research has often been epistemologically driven, guided by de Finetti's famous quote that ``probability does not exist''\cite{definetti1974}. For example, the ``quantum Bayesian'' methodology of Caves, Fuchs, \& Schack\cite{Caves2002a} has applied de Finetti's ideas to Bayes' theorem for use in quantum mechanics. In doing so, Caves et al. have argued that statistical systems are best interpreted by methods in which the Bayesian likelihood ratio is seen to be both external to the system and subjectively imposed on it by the observer\cite{Timpson2008}.\\

From a decision-making perspective the Caves et al.~approach is problematic. Bayes' theorem and, in particular, the ``naive Bayes' classifier'' have been used extensively to interpret information systems and develop normative decision-making models\cite{Oaksford2007}. While subjectivity may play a role in a descriptive model of human decision-making, its use in normative analysis could suggest the presence of a cognitive ``homonculus'' with the power to influence decision outcomes. Yet at a human scale, for instance, an observer's belief as to the chances of a fair coin landing either ``heads'' or ``tails'' has no known effect. Rather, within normative decision-making theory, the ``heads:tails'' likelihood ratio of 0.5:0.5 is only meaningful when considered as a property of the coin's own internal statistical system rather than as some ephemeral and arbitrary qualia.\\

However, there are axiomatic difficulties associated with Bayes' theorem, notably its reliance upon the use of marginal probabilities in the absence of structural statistical information, such as estimates of covariate overlap. This has led some researchers to attempt to reconceptualise psychology and decision-making theory using quantum mechanics\cite{Khrennikov2009,Busemeyer2012} - an approach with intuitive merit given that both disciplines apply statistical axioms to analyse and interpret probabilistic systems. Yet, despite progress made in this area, the lack of an orthodox Copenhagen-based theoretical counterpoint to Caves et al.~has impeded the development of new normative models. It is this knowledge gap which this article aims to fill.
\section{The limits of Bayes' theorem}
Bayes' theorem is used to calculate the conditional probability of a statement, or hypothesis, being true given that other information is also true. It is usually written as
\begin{equation}
\label{bayes}
P(H_{i}|D)=\frac{P(H_{i})P(D|H_{i})}{\sum\limits_{j}P(H_{j})P(D|H_{j})} \ .
\end{equation}
Here, $P(H_{i}|D)$ is the conditional probability of hypothesis $H_{i}$ being true given that the information $D$ is true; $P(D|H_{i})$ is the conditional probability of $D$ being true if $H_{i}$ is true; and $\sum\limits_{j}P(H_{j})P(D|H_{j})$ is the sum of the probabilities of all hypotheses multiplied by the conditional probability of $D$ being true for each hypothesis\cite{Oaksford2007}.
\begin{spacing}{1.5}
\begin{equation}
\label{figure1}
\begin{array}{lcc}
\hline 
& \mbox{Street A }(H_1) & \mbox{Street B }(H_2) \\  \hline
\mbox{Number of houses }(n) & 10 & 10\\ \hline
\mbox{\% blue front door }(D) & 0.8 & 0.7\\
\mbox{\% not blue front door }(\bar{D}) & 0.2 & 0.3\\ \hline
\end{array}
\end{equation}
\end{spacing}
To exemplify using the contingency information in \eqref{figure1}, if one wishes to calculate which of two streets is most likely to have a house with a blue front door, then using Bayes' theorem it is trivial to calculate that Street A is most likely with a likelihood ratio of approximately 0.53:0.47,
\begin{align}
\label{ex1}
P(H_{1}|D)&=\frac{0.5\times0.8}{(0.5\times0.8)+(0.5\times0.7)} = \frac{8}{15}&\approx 0.533 \ , 
\notag \\
P(H_{2}|D)&=1-P(H_{1}|D) = \frac{7}{15}&\approx 0.467 \ ,
\end{align}
where $P(H_i) = 10/(10+10) = 0.5$ for both $i=1,2$.\\

However, difficulties arise in the use of Bayes' theorem for the calculation of likelihood ratios where there are multiple non-exclusive data sets. For instance, if the information in \eqref{figure1} is expanded to include data about the number of houses with a garage \eqref{figure2} then the precise covariate overlap (i.e., $D_1 \cap D_2$) for each street becomes an unknown.
\begin{spacing}{1.5}
\begin{equation}
\label{figure2}
\begin{array}{lcc}
\hline 
&\mbox{Street A }(H_1)&\mbox{Street B }(H_2)\\  \hline
\mbox{Number of houses }(n)&10&10\\ \hline
\mbox{\% blue front door }(D_1)&0.8&0.7\\
\mbox{\% with garage }(D_2)&0.6&0.5\\ \hline
\end{array}
\end{equation}
\end{spacing}
All that may be shown is that, for each street, the number of houses with both features forms a range described by \eqref{min_max}, where $n(H_i)$ is the total number of houses in street $i$, $n(D_1|H_i)$ is the total number of houses in street $i$ with a blue front door, and $n(D_2|H_i)$ is the total number of houses in street $i$ with a garage,
\footnotesize
\begin{align}
\label{min_max}
\notag
& n(D_1 \cap D_2|H_{i}) \in
\\
& 
\left\{
\hspace{-0.2cm}
\begin{array}{l}
\Bigl[n(D_1|H_i) + n(D_2|H_i) - n(H_i) \ , \ldots, \min(n(D_1|H_i), n(D_2|H_i))\Bigr] \\
\qquad \qquad  \mbox{ if }  n(D_1|H_i)+n(D_2|H_i) > n(H_i) \ ,
\quad \mbox{or}
\\
\Bigl[0 \ ,\ldots, \min(n(D_1|H_i), n(D_2|H_i))\Bigr] \\ 
\qquad \qquad \mbox{ if } n(D_1|H_i)+n(D_2|H_i) \leq n(H_i) \ . 
\end{array}
\right.
\end{align}
\normalsize
Specifically for \eqref{figure2} these ranges equate to
\begin{align}
\label{min_max_actual}
n(D_1 \cap D_2|H_{1})&\in\{4,\:5,\:6\} \ , \notag \\
n(D_1 \cap D_2|H_{2})&\in\{2,\:3,\:4,\:5\} \ .
\end{align}
The simplest approach to resolving this problem is to naively ignore any intersection, or co-dependence, of the data and to directly multiply the marginal probabilities\cite[for example]{Doherty1979,Feeney2008,Daddario2012}. Hence, given \eqref{figure2}, the likelihood of Street A having the greatest number of houses with both a blue front door and a garage would be calculated as
\begin{align}
\label{pseudo_calc}
P(H_{1}|D_1 \cap D_2)&=\frac{0.5\times0.8\times0.6}{(0.5\times0.8\times0.6)+(0.5\times0.7\times0.5)} \notag \\
&\approx 0.578 \ .
\end{align}
Yet, because the data intersect, this probability value is only one of a number which may be reasonably calculated. Alternatives include calculating a likelihood ratio using the mean value $\mu$ of the frequency ranges for each hypothesis
\begin{align}
\label{Bayes2}
P(\mu [n(D_1 \cap D_2|H_1)])&=\frac{1}{10}\times\frac13(4+5+6) = 0.5\ ,\notag \\ 
P(\mu [n(D_1 \cap D_2|H_2)])&=\frac{1}{10}\times\frac14(2+3+4+5) = 0.35 \notag \\
&\Rightarrow \: P(H_1|\mu D_1 \cap D_2) \approx 0.588 \ ;
\end{align}
and taking the mean value of the probability range derived from the frequency range
\begin{align}
\label{Bayes3} 
&\min P(H_1|D_1 \cap D_2) = \frac{4}{4+5}
\ , \notag \\ \: 
&\max P(H_1|D_1 \cap D_2)=\frac{6}{6+2} \notag \\
&\Rightarrow \: \mu[ P(H_1|D_1 \cap D_2) ] \approx0.597 \ .
\end{align}
Given this multiplicity of probability values, it would seem that none of these methods may lay claim to normativity. This problem of covariate overlap has, of course, been previously addressed within statistical literature. For instance, the ``maximum likelihood'' approach, of Dempster, Laird, \& Rubin\cite{Dempster1977}, has demonstrated how an ``expectation-maximization'' algorithm may be used to derive appropriate covariate overlap measures. Indeed, the mathematical efficacy of this technique has been confirmed by Wu\cite{Wu1983}. However, from a psychological standpoint, it is difficult to see how such an iterative methodology could be employed in human decision-making, given mankind's reliance on heuristics rather than cognitively expensive exactitudes\cite{Kahneman2011}. Since there is also little evidence that the naive Bayes' classifier forms any part of the human decision-making process\cite{Doherty1979}, the theoretical advancement of the psychology of decision-making demands a mathematical approach in which covariate overlaps can be automatically, and directly, calculated from contingency data.

\section{A quantum mechanical proof of Bayes' theorem for independent data}

Previously unconsidered, the quantum mechanical von Neumann axioms would seem to offer the most promise in this regard, since the re-conceptualization of covariate data as a quantum entangled system allows for statistical analysis with few, non-arbitrary assumptions. Unfortunately there are many conceptual difficulties that can arise here. For instance, a Dirac notation representation of \eqref{figure2} as a standard quantum superposition is
\small
\begin{align}
\label{superpos}
\ket{\Psi}&=\frac{1}{\sqrt{N}}\biggl[ \alpha \bigg( \sqrt{\frac{1}{3}}\ket{4}_{H1}+\sqrt{\frac{1}{3}}\ket{5}_{H1}+\sqrt{\frac{1}{3}}\ket{6}_{H1} \bigg)\notag \\
&+\beta \bigg( \sqrt{\frac{1}{4}}\ket{2}_{H2}+\sqrt{\frac{1}{4}}\ket{3}_{H2}+\sqrt{\frac{1}{4}}\ket{4}_{H2}+\sqrt{\frac{1}{4}}\ket{5}_{H2} \bigg) \biggr] 
\ .
\end{align}
\normalsize
In this example, \eqref{superpos} cannot be solved since the possible values of $D_1 \cap D_2$ for each hypothesis \eqref{min_max_actual} have been described as equal chance outcomes within a general superposition of $H_1$ and $H_2$, with the unknown coefficients $\alpha$ and $\beta$ assuming the role of the classical Bayesian likelihood ratio.\\

The development of an alternative quantum mechanical description necessitates a return to the simplest form of Bayes' theorem using the case of exclusive populations $H_i$ and data sets $D$, $\bar{D}$, such as given in \eqref{figure1}. Here, the overall probability of $H_1$ may be simply calculated as
\begin{equation}
\label{Bayes4}
P(H_1)=\frac{n(H_1)}{n(H_1)+n(H_2)} \ .
\end{equation}
The a priori uncertainty in \eqref{figure1} may be expressed by constructing a wave function in which the four data points are encoded as a linear superposition
\begin{align}
\label{superpos2}
\ket{\Psi}=&\alpha_{1,1}\ket{H_1 \otimes D} + \alpha_{1,2}\ket{H_1 \otimes \bar{D}}\notag \\
&+\alpha_{2,1}\ket{H_2 \otimes D} + \alpha_{2,2}\ket{H_2 \otimes \bar{D}} \ .
\end{align}
Since there is no overlap between either $D$ and $\bar{D}$ or the populations $H_1$ and $H_2$, each datum automatically forms an eigenstate basis with the orthonormal conditions
\begin{align}
\label{ortho1}
\left<H_1 \otimes D|H_1 \otimes D \right> = \left<H_1 \otimes \bar{D}|H_1 \otimes \bar{D} \right>&= 1 \notag \\
\left<H_2 \otimes D|H_2 \otimes D \right> = \left<H_2 \otimes \bar{D}|H_2 \otimes \bar{D} \right>&= 1 \notag \\
\mbox{all other bra-kets } &= 0 \ ,
\end{align}
where the normalization of the wave function demands that
\begin{equation}
\label{ortho2}
\left<\Psi|\Psi\right>=1 \ ,
\end{equation}
so that the sum of the modulus squares of the coefficients $\alpha_{i,j}$ gives a total probability of 1
\begin{equation}
\label{ortho3}
|\alpha_{1,1}|^2 + |\alpha_{1,2}|^2 + |\alpha_{2,1}|^2 + |\alpha_{2,2}|^2 = 1 \ .
\end{equation}
For simplicity let
\begin{align}
\label{griddef}
x_1=P(D|H_1),\:\:y_1=P(\bar{D}|H_1) \ , \notag \\
x_2=P(D|H_2),\:\:y_2=P(\bar{D}|H_2) \ , \notag \\
X_1=P(H_1),\:\:X_2=P(H_2) \ .
\end{align}
If the coefficients $\alpha_{i,j}$ from \eqref{superpos2} are set as required by \eqref{figure1}, it follows that
\begin{equation}
\label{ortho4}
|\alpha_{1,1}|^2=x_1,\:\:|\alpha_{1,2}|^2=y_1,\:\:|\alpha_{2,1}|^2=x_2,\:\:|\alpha_{2,2}|^2 = y_2 \ ,
\end{equation}
so that the normalised wave function $\ket{\Psi}$ is described by
\begin{align}
\label{ortho6}
\ket{\Psi}=\frac{1}{\sqrt{N}}(&\sqrt{x_1}\ket{H_1\otimes D}+\sqrt{y_1}\ket{H_1\otimes\bar{D}}\notag \\
+&\sqrt{x_2}\ket{H_2\otimes D}+\sqrt{y_2}\ket{H_2\otimes\bar{D}})
\end{align}
for some normalization constant $N$.\\

The orthonormality condition \eqref{ortho2} implies that
\begin{equation}
\label{ortho7}
N=x_1+y_1+x_2+y_2=X_1+X_2 \ ,
\end{equation}
thereby giving the full wave function description
\footnotesize
\begin{equation}
\label{ortho8}
\ket{\Psi}=\frac{\sqrt{x_1}\ket{H_1\otimes D}+\sqrt{y_1}\ket{H_1\otimes\bar{D}}+\sqrt{x_2}\ket{H_2\otimes D}+\sqrt{y_2}\ket{H_2\otimes\bar{D}}}{\sqrt{X_1+X_2}} \ .
\end{equation}
\normalsize
If the value of $P(H_1|D)$ is to be calculated, i.e., the property $D$ is observed, then the normalized wave function \eqref{superpos2} necessarily collapses to
\begin{equation}
\label{collapse1}
\ket{\Psi'}=\alpha_1 \ket{H_1 \otimes D_1} + \alpha_2 \ket{H_2 \otimes D_1} \ ,
\end{equation}
where the coefficients $\alpha_{1,2}$ may be determined by projecting $\ket{\Psi}$ on to the two terms in $\ket{\Psi'}$ using \eqref{ortho1},  giving
\begin{align}
\label{collapse4}
&\alpha_1=\left<\Psi'|H_1\otimes D\right> = \sqrt{\frac{x_1}{X_1+X_2}} \ ,\notag \\
\mbox{\ } \notag \\
&\alpha_2=\left<\Psi'|H_2\otimes D\right> = \sqrt{\frac{x_2}{X_1+X_2}} \ .
\end{align}
Normalizing \eqref{collapse1} with the coefficient $N'$
\small
\begin{align}
\label{collapse2}
\ket{\Psi'}=\frac{1}{\sqrt{N'}}\Bigl(\sqrt{\frac{x_1}{X_1+X_2}} \ket{H_1 \otimes D} + \sqrt{\frac{x_2}{X_1+X_2}} \ket{H_2 \otimes D}\Bigr) \ ,
\end{align}
\normalsize
and using the normalization condition \eqref{ortho2}, implies that
\begin{align}
\label{ortho9}
1=\left<\Psi'|\Psi'\right>&=\frac{1}{N'}\Bigl( \frac{x_1}{X_1+X_2} + \frac{x_2}{X_1+X_2} \Bigr)\notag \\
\mbox{\ } \notag \\
\rightarrow N'&= \frac{x_1+x_2}{X_1+X_2} \ .
\end{align}
Thus, after collapse, the properly normalized wave function \eqref{collapse2} becomes
\begin{equation}
\label{collapse3}
\ket{\Psi'}=\sqrt{\frac{x_1}{x_1+x_2}} \ket{H_1 \otimes D} + \sqrt{\frac{x_2}{x_1+x_2}} \ket{H_2 \otimes D} \ ,
\end{equation}
which means that the probability of observing $\ket{H_1 \otimes D}$ is
\begin{align}
\label{bayesproof}
P(\ket{H_1\otimes D})=\Bigg(\sqrt{\frac{x_1}{x_1+x_2}}\:\Bigg)^2=\frac{\alpha_1^2}{\alpha_1^2+\alpha_2^2}=\frac{x_1}{x_1+x_2} \ .
\end{align}
This is entirely consistent with Bayes' theorem and demonstrates its derivation using quantum mechanical axioms.\\
\section{A quantum mechanical expression to calculate likelihood ratios, with co-dependent data}
Having established the principle of using a quantum mechanical approach for the calculation of simple likelihood ratios with mutually exclusive data \eqref{figure1}, it is now possible to consider the general case of $n$ hypotheses and $m$ data \eqref{figure3}, where the data are co-dependent, or intersect.
\begin{spacing}{1.5}
\begin{equation}
\label{figure3}
\begin{array}{|c|c|c|c|c|} \hline 
    & H_1 & H_2 & \cdots & H_n \\ \hline 
D_1 & x_{1,1} & x_{1,2} & \cdots & x_{1,n} \\ \hline
D_2 & x_{2,1} & x_{2,2} & \cdots & x_{2,n} \\ \hline
\vdots & & \vdots & & \vdots \\ \hline
D_m & x_{m,1} & x_{m,2} & \cdots & x_{m,n} \\ \hline
\end{array}
\end{equation}
\end{spacing}
Here the contingency table in \eqref{figure3} has been indexed using
\begin{equation}
x_{i, \alpha} \ , \qquad \alpha = 1, 2, \ldots, n; \quad i = 1,2, \ldots, m \ .
\end{equation}
While the general wave function remains the same as before, the overlapping data create non-orthonormal inner products which can be naturally defined as
\begin{equation}\label{innerGen}
\vev{H_\alpha \otimes D_i | H_\beta \otimes D_j} = 
c_{ij}^{\alpha} \delta_{\alpha\beta} \ , \:\:
c_{ij}^\alpha = c_{ji}^\alpha \in \IR \ , 
c_{ii}^\alpha = 1 \ ,
\end{equation}
Assuming, for simplicity, that the overlaps $c_{ij}^{\alpha}$ are real, then there is a symmetry in that $c_{ij}^\alpha = c_{ji}^\alpha$ for each $\alpha$. Further, for each $\alpha$ and $i$, the state is normalized, i.e., $c_{ii}^\alpha = 1$. The given independence of the hypotheses $H_\alpha$ also enforces the Kroenecker delta function, $\delta_{\alpha\beta}$.\\

The Hilbert space $V$ spanned by the kets $\ket{H_\alpha \otimes D_i}$ is $mn$-dimensional and, because of the independence of $H_\alpha$, naturally decomposes into the direct sum \eqref{directsum} with respect to the inner product, thereby demonstrating that the non-orthonormal conditions are the direct sum of $m$ vector spaces $V^\alpha$:
\begin{equation}
\label{directsum}
V = \mbox{Span}(\{ \ket{H_\alpha \otimes D_i} \} ) = 
\bigoplus\limits_{\alpha=1}^{n} V^{\alpha} \ , \qquad
\dim V^{\alpha} = m
\ .
\end{equation}
Since the inner products are non-orthonormal, each $V^{\alpha}$ must be individually orthonormalised. Given that $V$ splits into a direct sum, this may be achieved for each subspace $V^{\alpha}$ by applying the Gram-Schmidt algorithm\cite{Strang1980} to $\{ \ket{H_\alpha \otimes D_i} \}$ of $V$. Consequently, the orthonormal basis may be defined as
\begin{equation}\label{KHD}
\ket{K^\alpha_i} = \sum\limits_{k=1}^n A_{i,k}^{\alpha} \ket{H_\alpha \otimes D_k}
\ , 
\qquad
\vev{K^\alpha_i | K^\alpha_j} = \delta_{ij} \ ,
\end{equation}
for each $\alpha=1,2,\ldots,n$ with $m\times m$ matrices $A_{i,k}^{\alpha}$, for each $\alpha$.\\
Substituting the inner products \eqref{innerGen} gives
\begin{equation}\label{AAc}
\sum\limits_{k,k'=1}^m A_{ik}^{\alpha}  A_{jk'}^{\alpha} c_{kk'}^{\alpha}
= \delta_{ij} \quad \forall \alpha = 1, 2, \ldots, n \ .
\end{equation}
The wave-function may now be written as a linear combination of the orthonormalised kets $\ket{K^\alpha_i}$ with the coefficients $b_i^\alpha$, and may be expanded into the $\ket{H_\alpha \otimes D_i}$ basis using \eqref{KHD}, i.e.,
\begin{equation}\label{PsiGen}
\ket{\Psi} = \sum\limits_{\alpha,i} b_i^\alpha \ket{K^\alpha_i} =
\sum\limits_{\alpha,i,k} b_i^\alpha A_{ik}^{\alpha} \ket{H_\alpha \otimes D_k} \ .
\end{equation}
As with \eqref{ortho4} from earlier, the coefficients in \eqref{PsiGen} should  be set as required by the contingency table
\begin{equation}\label{bAx}
\sum\limits_{i} b_i^\alpha A_{i,k}^{\alpha} = \sqrt{x_{k \alpha}} \ ,
\end{equation}
where,  to solve for the $b$-coefficients, \eqref{AAc} may be used to invert
\begin{equation}
\sum\limits_{k,k'}\sum\limits_{i} b_i^\alpha A_{ik}^{\alpha} A_{jk'}c_{kk'}^\alpha
= \sum\limits_{k,k'} \sqrt{x_{k \alpha}} A_{jk'}^\alpha c_{k'k}^\alpha \ ,
\end{equation}
giving
\begin{equation}\label{bsolGen}
b_j^\alpha = \sum\limits_{k,k'} \sqrt{x_{k \alpha}} A_{jk'}^\alpha c_{kk'}^\alpha 
\ .
\end{equation}
Having relabelled the indices as necessary, a back-substitution of \eqref{bAx} into the expansion \eqref{PsiGen} gives
\begin{equation}
\ket{\Psi} = 
\sum\limits_{\alpha,i,k} 
b_i^\alpha A_{i,k}^{\alpha} \ket{H_\alpha \otimes D_k}  = 
\sum\limits_{\alpha,k} \sqrt{x_{k \alpha}} \ket{H_\alpha \otimes D_k} \ ,
\end{equation}
which is the same as having simply assigned each ket's coefficient to the square root of its associated entry in the contingency table.\\

The normalization factor for $\ket{\Psi}$ is simply $1/\sqrt{N}$, where $N$ is the sum of the squares of the coefficients $b$ of the orthonormalised bases $\ket{K^\alpha_i}$,
\begin{align}
\hspace{-0.5cm}
N &= \sum\limits_{i,\alpha} (b_i^\alpha)^2
= 
\sum\limits_{i,\alpha} 
b_i^\alpha
\left(
\sum\limits_{k,k'} \sqrt{x_{k \alpha}} A_{k',i}^\alpha c_{kk'}^\alpha 
\right) \notag \\
&=
\sum\limits_{k,k',\alpha} 
\sqrt{x_{k \alpha} x_{k'\alpha}} c_{kk'}^\alpha \ .
\end{align}
Thus, the final normalized wave function is
\begin{equation}
\label{NormalWave}
\ket{\Psi} =
\frac{\sum\limits_{\alpha,k} \sqrt{x_{k \alpha}} \ket{H_\alpha \otimes D_k}}
{\sqrt{\sum\limits_{i,j,\alpha} \sqrt{x_{i \alpha}x_{j \alpha}} c_{ij}^\alpha}}
\ ,
\end{equation}
where $\alpha$ is summed from 1 to $n$, and $i,j$ are summed from 1 to $m$.
Note that, in the denominator, the diagonal term  $\sqrt{x_{i \alpha}x_{j \alpha}}c_{ij}^\alpha$, which occurs whenever $i=j$, simplifies to $x_{i \alpha}$ since $c_{ii}^\alpha = 1$ for all $\alpha$. \\

From \eqref{NormalWave} it follows that, exactly in parallel to the non-intersecting case, if all properties $D_i$ are observed simultaneously, the probability of any hypothesis $H_\alpha$, for a fixed $\alpha$, is
\begin{align}\label{PHalpha}
 P(H_\alpha | D_1 \cap D_2 \ldots \cap D_m) &= \frac{\sum\limits_{i} (b_i^\alpha)^2}{\sum\limits_{i,\alpha} (b_i^\alpha)^2} 
=
\frac{
\sum\limits_{i,j} 
\sqrt{x_{i \alpha} x_{j\alpha}} c_{ij}^\alpha
}
{
\sum\limits_{i,j,\alpha} 
\sqrt{x_{i \alpha} x_{j\alpha}} c_{ij}^\alpha
}
\ ,
\end{align}
In the case of non-even populations for each hypothesis (i.e., non-even priors), each element within \eqref{PHalpha} should be appropriately weighted.
\section{Example solution}
Returning to the problem presented in the contingency table \eqref{figure2}, it is now possible to calculate the precise probability for a randomly selected house with the properties of ``blue front door'' and ``garage'' belonging to Street A ($H_1$). For this $2 \times 2$ matrix, recalling from \eqref{innerGen} that $c_{ii}^\alpha=1$ and $c_{ij}^\alpha = c_{ji}^\alpha$, the general expression \eqref{PHalpha} may be written as
\footnotesize
\begin{align}
& P(H_1|D_1 \cap D_2) 
= \frac{
\sum\limits_{i,j=1}^2 
\sqrt{x_{i,1} x_{j,1}} c_{ij}^1
}
{
\sum\limits_{i,j=1}^2  \sum\limits_{\alpha=1}^2  
\sqrt{x_{i \alpha} x_{j\alpha}} c_{ij}^\alpha
}
\notag \\
\mbox{\ } \notag \\
&= 
\frac{
\sqrt{x_{1,1}^2} c_{1,1}^1 + \sqrt{x_{2,1}^2} c_{2,2}^1 + 
\sqrt{x_{1,1}x_{2,1}} c_{1,2}^1 + \sqrt{x_{2,1}x_{1,1}} c_{2,1}^1 
}
{
\sum\limits_{\alpha=1}^2  
\sqrt{x_{1,\alpha}^2} c_{1,1}^1 + \sqrt{x_{2,\alpha}^2} c_{2,2}^1 + 
\sqrt{x_{1,\alpha}x_{2,\alpha}} c_{1,2}^1 + \sqrt{x_{2,\alpha}x_{1,\alpha}} c_{2,1}^1 
}
\notag \\
\mbox{\ } \notag \\
&=
\frac{x_1 + y_1 + 2 c_1 \sqrt{x_1 y_1}}{x_1 + x_2 + y_1 + y_2 + 2 c_1 \sqrt{x_1 y_1} + 2 c_2 \sqrt{x_2 y_2}} \ ,
\label{entang2}
\end{align}
\normalsize
where, adhering to the earlier notation \eqref{griddef},
\begin{align}
\label{griddef2}
x_1=x_{1,1} = P(D_1|H_1),\:\:y_1= x_{2,1} = P(D_2|H_1) \ , \notag \\
x_2=x_{1,2} = P(D_1|H_2),\:\:y_2= x_{2,2} = P(D_2|H_2) \ , \notag \\
X_1=P(H_1),\:\:X_2=P(H_2) \ ,\
\end{align}
and, for brevity, $c_1 := c_{1,2}^1$ , $c_2 := c_{1,2}^2$ .
For simplicity, $P(H_i|D_1 \cap D_2)$ will henceforth be denoted as $P_i$.
Implementing \eqref{entang2} is dependent upon deriving solutions for the yet unknown expressions $c_{i}$, $i=1,2$ which govern the extent of the intersection in \eqref{innerGen}. This can only be achieved by imposing reasonable constraints upon $c_i$ which have been inferred from expected behaviour and known outcomes, i.e., through the use of boundary values and symmetries. Specifically, these constraints are:
\begin{description}
\item[Data dependence.] The expressions $c_i$ must, in some way, be dependent upon the data given in the contingency table, i.e.,
\begin{align}
c_1&=c_1(x_1,y_1,x_2,y_2;X_1,X_2) \ , \notag \\
c_2&=c_2(x_1,y_1,x_2,y_2;X_1,X_2) \ .
\end{align}
\item[Probability.] The calculated values for $P_i$ must fall between 0 and 1. 
Since $x_i$ and $y_i$ are positive, it suffices to take
\begin{equation}
-1\:<\:c_i(x_1,\:y_1,\:x_2,\:y_2)\:<\:1
\ .
\end{equation}
\item[Complementarity.] The law of total probability dictates that
\begin{equation}
P_1+P_2=1 \ ,
\end{equation}
which can easily be seen to hold.
\item[Symmetry.] The exchanging of rows within the contingency tables should not affect the calculation of $P_i$. In other words, for each $i=1,2$, $P_i$ is invariant under $x_i \leftrightarrow y_i$. 
This constraint implies that
\begin{equation}
c_i(x_1,y_1,x_2,y_2)=c_i(y_1,x_1,y_2,x_2) \ .
\end{equation}
Equally, if the columns are exchanged then $P_i$ must map to each other, i.e., for each $i=1,2$ then $P_1 \leftrightarrow P_2$ under $x_1 \leftrightarrow x_2,\: y_1 \leftrightarrow y_2 $ which gives the further constraint that
\begin{equation}
c_1(x_1,y_1,x_2,y_2)=c_2(x_2,y_2,x_1,y_1) \ .
\end{equation}
\item[Known values.] There are a number of contingency table structures which give rise to a known probability, i.e.,
\begin{spacing}{1.5}
\begin{align}
\begin{array}{|c|c|l|} \hline 
    & H_1 & H_2  \\ \hline
D_1 & 1 & 1  \\
D_2 & m & n  \\ \hline
\end{array} \notag
&\ \ \ \ \ \ \ \ \ \rightarrow
&P_1 &= \frac{m}{m+n}
\\ \notag
\begin{array}{|c|c|l|} \hline 
    & H_1 & H_2  \\ \hline
D_1 & m & n  \\
D_2 & 1 & 1  \\ \hline
\end{array} \notag
&\ \ \ \ \ \ \ \ \ \rightarrow
&P_1 &= \frac{m}{m+n}
\\
\begin{array}{|c|c|l|} \hline 
    & H_1 & H_2  \\ \hline
D_1 & n & m  \\
D_2 & m & n  \\ \hline
\end{array} \notag
&\ \ \ \ \ \ \ \ \ \rightarrow
&P_1 &= \frac{1}{2}
\\
\begin{array}{|c|c|l|} \hline 
    & H_1 & H_2  \\ \hline
D_1 & n & n  \\
D_2 & m & m  \\ \hline
\end{array} \notag
&\ \ \ \ \ \ \ \ \ \rightarrow
&P_1 &= \frac{1}{2}
\\
\begin{array}{|c|c|l|} \hline 
    & H_1 & H_2  \\ \hline
D_1 & m & m  \\
D_2 & m & m  \\ \hline
\end{array}
&\ \ \ \ \ \ \ \ \ \rightarrow
&P_1 &= \frac{1}{2} \ ,
\end{align}
\end{spacing}
where $m,n$ are positively valued probabilities. For such contingency tables the correct probabilities should always be returned by $c_i$. Applying this principle to \eqref{entang2} gives the constraints
\footnotesize
\begin{align}
\label{movern}
&\frac{m}{m+n} = \frac{2 c_1(m,1,n,1) \sqrt{m}+m+1}
   {2 c_1(m,1,n,1) \sqrt{m}+2 c_2(m,1,n,1)\sqrt{n}+m+n+2} \ , \\ \notag
\mbox{\ } \notag \\
\label{movern2}
&\frac{1}{2} = \frac{2 c_1(n,m,m,n) \sqrt{m} \sqrt{n}+m+n}{2 c_1(n,m,m,n) \sqrt{m} \sqrt{n}+2 c_2(n,m,m,n) \sqrt{m} \sqrt{n}+2 m+2 n} \ ,
\\ \notag
\mbox{\ } \notag \\
\label{movern3}
&\frac{1}{2} = \frac{2 c_1(n,m,n,m) \sqrt{m} \sqrt{n}+m+n}{2 c_1(n,m,n,m) \sqrt{m} \sqrt{n}+2 c_2(n,m,n,m) \sqrt{m} \sqrt{n}+2 m+2 n} \ .
\end{align}
\normalsize
\item[Non-homogeneity.] Bayes' theorem returns the same probability for any linearly scaled contingency tables, e.g.,
\footnotesize
\begin{align}
\label{homogen1}
&x_1 \rightarrow 1.0,\: y_1 \rightarrow 1.0,\: x_2 \rightarrow 1.0,\: y_2 \rightarrow 0.50 \Rightarrow P_1\approx 0.667 \ ,\\
\label{homogen2}
&x_1 \rightarrow 0.5,\: y_1 \rightarrow 0.5,\: x_2 \rightarrow 0.5,\: y_2 \rightarrow 0.25 \Rightarrow P_1\approx 0.667 \ .
\end{align}
\normalsize
While homogeneity may be justified for conditionally independent data, this is not the case for intersecting, co-dependent data since the act of scaling changes the nature of the intersections and the relationship between them. This may be easily demonstrated by taking the possible value ranges for \eqref{homogen1} and \eqref{homogen2}, calculated using \eqref{min_max}, which are
\begin{align}
\mbox{Eq.~\eqref{homogen1}} \Rightarrow &(D_1 \cap D_2)|H_{1}=\{1\} \ , \notag \\
&(D_1 \cap D_2)|H_{2}=\{0.5\} \ , \notag \\
\mbox{\ } \notag \\
\mbox{Eq.~\eqref{homogen2}} \Rightarrow &(D_1 \cap D_2)|H_{1}=\{0.0 \ldots 0.5\} \ , \notag \\
&(D_1 \cap D_2)|H_{2}=\{0.0 \ldots 0.25\} \ .
\end{align}
The effect of scaling has not only introduced uncertainty where previously there had been none, but has also introduced the possibility of 0 as a valid answer for both hypotheses. Further, the spatial distance between the hypotheses has also decreased. For these reasons it would seem unreasonable to assert that \eqref{homogen1} and \eqref{homogen2} share the same likelihood ratio.
\end{description}
Using these principles and constraints it becomes possible to solve $c_i$. From the principle of symmetry it follows that
\begin{align}
c_1(n, m, m, n) = c_2(m,n,n,m) = c_2(n, m, m, n) \ , \notag
\\
c_1(n, m, n, m) = c_2(n,m,n,m) = c_2(n, m, n, m) \ ,
\end{align}
and that the equalities \eqref{movern2}, \eqref{movern3} for $P_i=0.5$ automatically hold. Further, \eqref{movern} solves to give
\begin{equation}
\label{function2}
c_2(m,1,n,1) = \frac{2 \sqrt{m} n c_1(m,1,n,1) - m + n}{2m \sqrt{n}} \ ,
\end{equation}
which, because $c_1(n,1,m,1)=c_2(m,1,n,1)$, finally gives
\begin{equation}
\label{function1}
c_1(n,1,m,1) =
\frac{2 \sqrt{m} n c_1(m,1,n,1) - m + n}{2m \sqrt{n}} \ .
\end{equation}
Substituting $g(m,n) := \sqrt{n} c_1(m,1,n,1)$
transforms \eqref{function1} into an anti-symmetric bivariate functional equation in $m,n$,
\begin{align}
g(m,n) - g(n,m) = \frac{m}{2\sqrt{mn}} - \frac{n}{2\sqrt{mn}} \ ,
\end{align}
whose solution is $g(m,n)=\frac{m}{2\sqrt{mn}}$\ . \\

This gives a final solution for the coefficients $c_{1,2}$ of
\begin{align}\label{c12sol}
c_1(x_1,y_1,x_2,y_2)&=\frac{\sqrt{x_1 y_1}}{2x_2 y_2} \ , \notag \\ \quad
c_2(x_1,y_1,x_2,y_2)&=\frac{\sqrt{x_2 y_2}}{2x_1 y_1} \ .
\end{align}
Thus, substituting \eqref{c12sol} into \eqref{entang2} gives the likelihood ratio expression of,
\begin{equation}
P(H_1|D_1 \cap D_2) =\frac{\frac{x_1 y_1}{x_2 y_2}+x_1+y_1}
{
\frac{x_1 y_1}{x_2 y_2}+ x_1 + y_1 +
\frac{x_2 y_2}{x_1 y_1}+ x_2 + y_2
} \ .
\end{equation}
Given that the population sizes of $H_1$ and $H_2$ are the same, no weighting of the elements needs to take place. Hence, the value of $P(H_1|D_1 \cap D_2)$ for \eqref{figure2} may now be calculated to be
\begin{equation}
P(H_1|D_1 \cap D_2) \approx 0.5896 \ .
\end{equation}
\section{Discussion}
One of the greatest obstacles in developing any statistical approach is demonstrating correctness. This formula is no different in that respect. If correctness could be demonstrated then, a priori, there would be an appropriate existing method which would negate the need for a new one. All that may be hoped for in any approach is that it generates appropriate answers when they are known, reasonable answers for all other cases, and that these answers follow logically from the underlying mathematics.\\

However, what is clear is that the limitations of the naive Bayes' classifier render any calculations derived from it open to an unknown margin of error. Given the importance of accurately deriving likelihood ratios this is troubling. This is especially true when these calculations are used to describe normative psychological theories \cite{Doherty1979} from which inferences are drawn as to the failure of human logic.\\

As a quantum mechanical methodology this result is able to calculate accurate, iteration free, likelihood ratios which fall beyond the scope of existing statistical techniques, and offers a new theoretical approach within cognitive psychology. Further, through the addition of a Hamiltonian operator to introduce time-evolution, it can offer likelihood ratios for future system states with appropriate updating of the contingency table. In contrast, Bayes' theorem is unable to distinguish directly between time-dependent and time-independent systems. This may lead to situations where the process of contingency table updating results in the same decisions being made repeatedly with the appearance of an ever increasing degree of certainty. Indeed, from \eqref{bayesproof}, it would seem that the naive Bayes' classifier is only a special case of a more complex quantum mechanical framework, and may only be used where the exclusivity of data is guaranteed.\\

The introduction of a Hamiltonian operator, and a full quantum dynamical formalism, is in progress, and should have profound implications for cognitive psychology. Inevitably, such a formalism will require a sensible continuous classical limit. In other words, the final expressions for the likelihood ratios should contain a parameter, in some form of $\hbar$, which, when going to 0, reproduces a classically known result. For example, the solutions to \eqref{c12sol} could be moderated as
\begin{align}
c_1(x_1,y_1,x_2,y_2) = \frac{\sqrt{x_1y_1}}{2x_2y_2}(1 - \exp(-\hbar)) \ , \notag \\
\mbox{\ } \notag \\
c_2(x_1,y_1,x_2,y_2) = \frac{\sqrt{x_2y_2}}{2x_1y_1}(1 - \exp(-\hbar)) \ ,
\end{align}
so that in the limit of $\hbar \rightarrow 0$, the intersection parameters, $c_1$ and $c_2$, vanish to return the formalism to the classical situation of independent data.
\section{Conclusion}
This article has demonstrated both theoretically, and practically, that a quantum mechanical methodology can overcome the axiomatic limitations of classical statistics in respect to their application within cognitive psychology. In doing so, it challenges the orthodoxy of de Finetti's epistemological approach to statistics by demonstrating that it is possible to derive ``real'' likelihood ratios from information systems without recourse to arbitrary and subjective evaluations.\\

While further theoretical development work needs to be undertaken, particularly with regards to the application of these mathematics in other domains, it is hoped that this article will help advance the debate over the nature and meaning of statistics, particularly with respect to cognitive psychology and decision-making theory.

\subsection{\label{sec:citeref}Citations and References}
\bibliographystyle{apsrev4-1}
\bibliography{Bond_He_Ormerod_2015_MS_PRX}

\end{document}